\documentclass[10pt]{amsart}
\usepackage{type1ec}
\usepackage{latexsym}
\usepackage{graphicx}
\usepackage{verbatim}

\input epsf

\begin{document}

\medskip
\title
[The Mathematical Background of Lomonosov's Contribution]
{
The Mathematical Background of\\ Lomonosov's Contribution
}
\author{S.~S. Kutateladze}
\address[]{
Sobolev Institute of Mathematics\newline
\indent 4 Koptyug Avenue\newline
\indent Novosibirsk, 630090\newline
\indent RUSSIA
}
\begin{abstract}
This is a short overview of the influence of mathematicians and
their ideas on the
creative contribution of Mikha\u\i{}lo Lomonosov on the occasion
 of the tercentenary of his birth.
\end{abstract}
\email{
sskut@math.nsc.ru
}
\keywords{Euler, Leibniz, Lomonosov, Newton,  Wolff,
 monad, calculus, nonstandard analysis}

\date{April 12, 2011}
\maketitle

%\baselineskip=.95\baselineskip

Lomonosov is the Russian colossus  of the epoch of scientific giants.
Lomonosov was not a mathematician, but without mathematicians there would
be no Lomono\-sov as the first and foremost Russian scientist at all.

Science in Russia had started with the foundation of the Academy of Sciences and Arts
which then evolved  into the Russian Academy of Sciences of these days.
The turn of the sixteenth and seventeenth centuries is  a signpost of the
history of the mankind, the onset of the organized science.
The time of the birth  of scientific societies and academies
accompanied the revolution in the natural sciences which
rested upon the discovery of differential and integral calculus.
The new language of mathematics brought about
an opportunity to make impeccably precise predictions of
future events.

To the patriotism of Peter the Great and the cosmopolitanism of Leibniz
we owe the foundation of the Saint Petersburg Academy of Sciences as
the center of Russian science. Peter and Leibniz stood at the cradle
of Russian science in much the same way as Euler and Catherine~I
are the persons from whom we count the history of the national mathematical school
in Russia.  We  must also acclaim the  outstanding role
of Leibniz who  prepared for Peter  a detailed plan of organizing  academies
in Russia  (cp.~\cite{Anri}).
Leibniz viewed Russia as a bridge for connecting Europe with China whose
Confucianism  would inoculate some necessary  ethical principles for
bringing moral health to Europe (cp.~\cite{Hellenbroich}).
Peter wanted to see Leibniz as an active organizer of the Saint Petersburg Academy,
he persuaded Leibniz in person and appointed Leibniz a {\it Justizrat} with a lavish salary.
Elisabeth Charlotte d'Orl\'eans who was present at the meeting of Peter with Leibniz
wrote on December 10, 1712 (cp.~\cite{Smith}):

\begin{itemize}
\item[]{\small
Muscovy must be a savage place. Thus I find that Herr Leibniz
is right in not wishing to go there. I am as it were charmed by
the Tsar, when I see how much care he takes to improve his country.}
\end{itemize}

 It is worth observing that Peter visited   the Royal Mint in London
in~1698 during the so-called  ``Grand Embassy.''
At that time Newton was  Warden of the Mint and we can hardly imagine that
he ignored Peter's visits. Nevertheless there is no evidence that Peter met
Newton. It is certain that Jacob Bruce, one of the closest  associates  of Peter,
had discussions with Newton (cp.~\cite[p.~199]{Cracraft}).
In 1714, two years after Peter made Leibniz a  {\it Justizrat},
Aleksandr Menshikov applied for membership in the London Royal Society, which was
an extraordinary and unforeseeable event. What is more mysterious, Menshikov's
application was approved and he was notified of his new status by a letter from Newton
himself (cp.~\cite[Ch.~16]{Vavilov}).
	
The genius of Newton has revealed to the universe the mathematical
laws of nature, disclosed to a mathematician
a universal language for describing the ever-changing world.
The genius of Leibniz has pointed out to the mankind  the opportunities
of mathematics as a reliable method of reasoning, the genuine logic
of human knowledge.
Leibnizian ideas of {\it mathesis universalis\/}  and  {\it calculemus}
arose once and forever as a dream and instrument of
science.

The influence of the ideas of Newton and Leibniz resulted in the
the scientific outlook (e.g.,~\cite[Ch.~2]{Ekeland}). The revolt of the
natural sciences at the turn of the seventeenth and eighteenth centuries
was determined  by the invention of differential and integral calculus.
The competing ideas of the common mathematics of Newton and Leibniz
determined all principal trails of thought  of the intellectual search of the epoch.
The contribution of   Lomonosov exhibits a convincing example of the general trends.
To grasp the scientific approaches of Lomonosov,  to understand his creative
revelations and naive delusions is impossible without deep analysis of
and thorough comparison between the views of Newton and Leibniz.

The monads of Leibniz as well as the fluxions and fluents of Newton
are  products of the heroic epoch of the telescope and microscope.
The independence of the discoveries of Leibniz and Newton  is obvious, since
their approaches, intellectual backgrounds, and intentions
were radically different. Nevertheless, the groundless priority quarrel between
Leibniz and Newton  has become the behavioral pattern
for many generations of scientists.
Leibniz and Newton discovered the same formulas, part of which had already been known.
Leibniz, as well as Newton, had his own priority in the invention of
differential and integral calculus. Indeed, these scientists
suggested the versions of mathematical analysis which were based
on different grounds. Leibniz founded analysis on actual infinitesimals,
erecting the tower of his perfect philosophical system known as monadology. The key of Newton was
his method of ``prime and ultimate ratios'' which is rightfully
associated with the modern limit theory.

The Leibnizian stationary vision of mathematical objects counterpoises the
Newtonian dynamical perception of ever-changing variable quantities.
The source of the ideas by Leibniz was the geometrical views
of antiquity which he was enchanted with from his earliest infantry.
The monad of Euclid is the mathematical tool of calculus, presenting a twin
to the point, the atom of geometry. Mathematics of Euclid is the product
of the human spirit. The monads of Leibniz, nurtured by his dream of {\it calculemus}
are the universal instrument of creation whose understanding
brings a man to the divine providence in creating the best of all possible worlds.

The point and the monad of the ancients  are independent forms of reasoning,
mental reflections of indivisible constituents of figures and numbers.
Both ideas are tightly woven into the conception of universal atomism.
The basic idea of the straight line  has  incorporated
the understanding of its dualistic---discrete-continuous---nature from the very
beginning of geometry. Leibniz ascribed the universal meaning to the ancient geometrical idea,
discerning the divine providence  that is incorporated in the idea.

Newton  got acquaintance with Euclid only in his ripe years and so he travelled in his own way,
perceiving universal motion as something done at the creation of the world that
could thus never be reduced to any sum of states of rest.
The perfectly precise characterization of Newton was done by Keynes
in his talk~\cite{Keynes} prepared to the tercentenary of Newton which had
to be celebrated in 1942, but was postponed until 1946 because of
the circumstances of wartime. Unfortunately, Keynes had passed away
three months before the celebration and his lecture was delivered by his brother.
Keynes wrote:

\begin{itemize}
\item[]{\small
Why do I call him a magician? Because he looked on the whole universe and
all that is in it as a riddle, as a secret which could be read by applying
pure thought to certain evidence, certain mystic clues which God had laid
about the world to allow a sort of philosopher's treasure hunt to the esoteric
brotherhood. He believed that these clues were to be found partly in the evidence
of the heavens and in the constitution of elements (and that is what gives
the false suggestion of his being an experimental natural philosopher), but also partly
in certain papers and traditions handed down by the brethren in an unbroken chain back
to the original cryptic revelation in Babylonia. He regarded the universe as a cryptogram
set by the Almighty---just as he himself wrapt the discovery of the calculus in a cryptogram
when he communicated with Leibniz. By pure thought, by concentration of mind, the riddle,
he believed, would be revealed to the initiate.''
}
\end{itemize}

If Newton was the last scientific magician, then Leibniz was
the first mathematical dreamer.

The outlook of Leibniz, proliferating with his works,
occupies a unique place in  human culture.
We can hardly find in the philosophical treatises of his
predecessors and later thinkers something
comparable with the phantasmagoric conceptions of monads,
the special and stunning
constructs of the world and mind which precede, comprise, and
incorporate all the infinite advents of the eternity.
It is worth emphasizing that mathematics was the true source of
of the philosophical ideas of Leibniz. Suffice it to quote  Child  who translated into English
and commented the early mathematical papers of Leibniz (cp.~\cite[Preface]{Child}):

\begin{itemize}
\item[]{\small
The main ideas of his philosophy are to be attributed to his mathematical
work, and not {\it vice versa}.
}
\end{itemize}

{\it Monadology} \cite[pp.~413--428]{Mon} is usually dated as of 1714.
This article was never published during Leibniz's life.
Moreover, it is generally accepted  that the very term ``monad''
had appeared in his writings since 1690 when he was already an
established and prominent scholar.

The special attention to the origin of the
term ``monad'' and the particular investigation into
the date of its first appearance in the works by Leibniz are in
fact the present-day products.
There are now a few if any cultivated persons who never
got acquaintance with the basics of planimetry and heard nothing of
Euclid. However, no one has ever met the concept of ``monad''
on the school bench.
Neither the contemporary translations of Euclid's {\it Elements}
nor the popular school text-books
contain this seemingly exotic term. However,
the concept of  ``monad''  is fundamental not only
for Euclidean geometry but also for the whole science of
the Ancient Hellas.

By Definition I of Book VII of Euclid's {\it Elements} \cite{Euc}
a~monad is ``that by virtue of which each
of the things that exist is called one.''
Euclid proceeds with Definition~2:
``A number is a multitude composed of monads.''
Note that the present-day translations of the Euclid treatise
substitute ``unit'' for ``monad.''

A contemporary reader can hardly understand why
Sextus Empiricus, an outstanding scepticist
of the second century, wrote when
presenting the mathematical views of his predecessors
as follows \cite{Sextus}:
``Py\-tha\-goras said that the origin of the
things that exist is a~monad by virtue of which
each of the things that exist is called one.''
And furthermore:
``A~point is structured as a~monad; indeed,  a~monad is
a~certain origin of numbers and likewise a~point is a~certain origin
of lines.''
Now some place is in order for the excerpt which
can easily be misconceived as a citation
from {\it Monadology}:
``A~whole as such is indivisible and
a~monad, since it is a~monad, is not divisible.
Or, if it splits into many pieces
it becomes a~union of many monads
rather than a~[simple] monad.''

It is worth observing that the ancients
sharply perceived an exceptional status of the
start of counting. In order to count, one
should firstly  particularize the entities to count and
only then to proceed with putting these entities into correspondence
with some symbolic series of numerals.
We begin counting with making ``each of the things one.''
The especial role of the start of counting
is reflected in the almost millennium-long dispute
about whether or not the unit (read, monad) is a natural number.
We feel today that it is excessive to distinguish the key role
of the unit or monad which signifies the start of counting.
However, this was not always so.

From the times of Euclid, all serious scientists knew
about existence of the two basic concepts of mathematics: a point and a monad.
By Definition~1 of Book~1 of Euclid's {\it Elements}:
``A point is that which has no parts.''
Clearly this definition differs drastically from the definition of
monad as that which makes one from many.
The cornerstone  of geometry is other than that of arithmetic.
Without clear understanding of this circumstance it is impossible
to comprehend the essence of the views of Leibniz.
By the way, the modern set theory refers to ``that which has no parts''
as the empty set, the starting cardinal of the von Neumann universe.
The present-day mathematics seems to have no concept
that is vocalized as ``that which many makes into one.''
We will return to the modern mathematical definition of monad
shortly.

As a top mathematician of his epoch, Leibniz
was in full command of Euclidean geometry.
Therefore, rather bewildering is  Item~1 of  {\it Monadolody}
where Leibniz gave the first idea of what his monad actually is:

\begin{itemize}
\item[]{\small
The Monad, of which we shall here speak, is nothing but a
simple substance, which enters into compounds. By ``simple'' is meant
``without parts.''
}
\end{itemize}

This definition of monad as a ``simple'' substance
without parts coincides with the Euclidean definition of
point. At the same time the reference to compounds consisting of monads
reminds us the structure of the definition of number which belongs to Euclid.

The synthesis of both primary definitions of Euclid
in the Leibnizian monad is not accidental.
We must always bear in mind that the seventeenth century
is the epoch of microscope. It was already in the 1610s that microscopes
were mass-produced in many European countries. From the 1660s
Europe was enchanted by Antony van Leeuwenhoek's microscope.

Attempting to pursue the way of Leibniz's thought,
we must always keep in mind that he was a mathematician by
belief. From his earliest childhood,  Leibniz dreamed of
``some sort of calculus'' that operates in the
``alphabet of human thoughts'' and possesses
the same beauty, strength, and integrity as mathematics in solving
arithmetical and geometrical problems. Leibniz devoted many articles
to invention of this universal logical calculus.
He remarked that his general methodological views
are grounded on the ``studying of the ways of analysis in mathematics
to which I was subjected with such an ardency that I do not know whether
there are many to be found today who invested much more toil into  it than me.''

The teacher of Lomonosov was Christian Wolff,  an ardent propagator of the ideas of monadism
and the mathematical method.   Wolff was considered by his contemporaries as
the second figure after Leibniz in the continental science.
The first figure of  Misty Albion was Newton.
It is impossible to forget that the intellectual life of that epoch was
heavily contaminated with the nasty controversy about priority between
Newton and Leibniz.
The deplorable consequence of the confrontation was the stagnation
and isolation of the mathematical life of England.
As regards the continent, the slight but perceptible  neglect
to the contribution by Newton led to dogmatization and canonization
of the teaching of Leibniz which was often understood with distortions.

Wolff was an epigone rather than a follower of Leibniz.
Tore Fr\"angsmyr remarked in~\cite[p.~34]{Spirit}:

\begin{itemize}
\item[]{\small
Wolff's epistemology is quite simple. Human knowledge outside
Christian revelation can be acquired in three ways: by experience
(historical knowledge), by reason (mathematical knowledge), or by a
combination of the two (philosophical knowledge). The last of these
three methods is preferable; the other two have value only to the
extent that they can be of use to philosophy.}
\end{itemize}

The true disciples  of Leibniz  were Jacob, Jean, and  Jacques Bernoulli as well
as Euler who was a self-taught prodigy close to Bernoulli by the vogue and
understanding of life.

Observe that  Wolff was the trendsetter in
mathematical education of the beginning of the eighteenth century.
After Leibniz's refusal to transfer to Saint Petersburg
for organizing the Academy, Peter considered Wolff as its possible leader.
Wolff's  treatise {\it Der Anfangsgr\"unde aller mathematischen Wissenshaft\/}
was published in four parts in~1710, abridged later for a wider readership,
and reprinted many times~(cp.~\cite[p.~23]{Yushkevich}).

Explicating his pedagogical principles, Wolff wrote
(cp.~\cite{Rogers}):\footnote{Unfortunately, the next quotation is
not a direct translation of Wolff's original.}

\begin{itemize}
\item[]{\small
In my lectures I paid most attention to the three aspects:\\
 (1) I never used any word that was left unexplained before in order
 to avoid ambiguity and logical gaps;\\
(2) I never used any theorem I had not proved before;\\
(3) I always connected theorems and definitions with one another to make
a continual logical chain.}
\item[]{\small
It is universally known that these rules are
followed in mathematics. If we compare the mathematical
method of education  with the approach of logic, we can see that
the mathematical method of education is nothing other than
the exact application of the rules of inference.
Therefore it is immaterial whether we use the mathematical method
of education or the rules of inference provided that the latter
are true. Inasmuch I have demonstrated that the mathematical reasoning
reflects the natural reasoning, and the logical reasoning
is just a definitely improved form of the natural reasoning,
I have any right to declare that my method of education follows
the natural mode of reasoning.}
\end{itemize}

Hegel was rather sceptical about the pedagogical style of Wolff and
remarked (cp.~\cite[p.~363]{Hegel}):

\begin{itemize}
\item[]{\small
 Wolff on the one hand started upon a large range of investigation,
 and one quite indefinite in character, and on the other, held to
 a strictly methodical manner with regard to propositions and their proofs.
 The method is really similar to that of Spinoza, only it is more wooden
 and lifeless than his.}
\end{itemize}

The ideas of Wolff in education were
well accepted by Lomonosov, since Wolff and he were connected
with the warm feeling of mutual respect.
Wolff's  mathematical method was a basis of
Lomonosov's scientific articles during many years of his creativity.
It should be observed that , unlike Wolff who had excellent
mathematical training, Lomonosov was not sufficiently well acquainted
even with Euclid's {\it Elements}, and he never possessed a working knowledge
of differential and integral calculus.

We must emphasize that  Lomonosov never met
Euler (mentioning this in his famous talk ``Lomonosov and
World Science'' \cite{Kapitsa}, Kapitsa made an exquisite circumlocution
``of course we cannot exclude the possibility of Lomonosov presence at
the public lecture of Euler that he delivered before his departure to Germany'').
These circumstances explain to us  we can hardly  find any practical
applications of mathematics in the papers of Lomonosov and why
some of his thoughts about the nature of mathematical knowledge
are naive and incorrect.

For instance, in his great paper  
{\it Meditationes
de Caloris et Frigoris Causa
Auctore Michaele Lomonosow} propounding foundations of the
molecular-kinetic  theory of heat~(cp.~\cite[Cp.~1]{KTS}),
Lomonosov wrote~\cite[c.~24]{Tom2}:

\begin{itemize}
\item[]{\small
Nulla demonstrandi methodus certior est ea mathematicorum,
qui deductas a priori propositiones exemplis vel examine
institute a posteriore confirmare solent.}
\vskip2pt
\item[]{\small
There is no more certain method of demonstration than
the method of mathematicians who corroborate
the  {\it a priori} deduced
propositions by example and {\it a posteriori} examination.}
\end{itemize}

It is worth emphasizing that from this formally wrong thesis about
the nature of mathematical proof, Lomonosov deduced
the remarkable and undoubtedly true conclusion:

\begin{itemize}
\item[]{\small
Idcirco nostram theoriam
ulterius prosecuturi, ad exemplum eorum phaenomena praecipua,
quae circa ignem et calorem observantur, explicando
assertum \S\ 11 verissimum esse confirmabimus.}
\vskip2pt
\item[]{\small
Therefore, to further develop our theory we use the lead of mathematicians
and explain the most important phenomena that are observed for fire and heat,
thus corroborating  the full validity of the assertion we made in
\S 11.}
\end{itemize}

In actuality, Lomonosov discussed in this excerpt  the technology of mathematical modeling
which differs drastically from any mathematical formalism as such.

The attitude of Lomonosov to monads deserves a slightly thorough examination.
Developing the atomistic ideas of corpuscular physics in his papers
of 1743 and 1744, i.~e. {\it Tentamen Theoriae
De Particulis Insensibiubus Corporum
Deque Causis Qualitatum Particularium in Genere},
{\it De Cohaesione Et Situ Monadum Physicarum}, and
{it De Particulis Physicis Insensibilibus
Corpora Naturalia Constituentibus,
in Quibus Qualitatum Particularium
Ratio Suffic1ens Continetur} (\cite[pp. 169--235, 265--314]{Tom1})
as well as in his extensive correspondence,  Lomonosov
sparingly use the concept of monad, especially distinguishing
{\it  monades physicae}.
The physical monads of Lomonosov are closer to the conception of atoms
rather than mathematical monads or Leibnizian ideal monads.
Long-term personal contemplations over the structure of the matter led Lomonosov 
(e.g., see \cite{Gareth Jones}). This is reflected in the choice of
Latin scientific terms of the later papers of Lomonosov~(cp.~\cite{Terms}).

In February of 1754 Lomonosov wrote to Euler \cite[pp.~501-502]{Tom10}:

\begin{itemize}
\item[]{\small
Fateor me idcirco potius ilia praeteriisse,
ne magnorum virorum scripta invadendo sui ostentator viderer
potius, quam veritatis scrutator. Haec ipsa ratio jam longo
tempore prohibet, quominus meditationes meas de monadibus
erudito orbi proponam discutiendas. Quamvis enim mysticam
fere illam doctrinam funditus everti argumentis meis debere
confidam; viri tamen, cujus erga me officia oblivisci non
possum, senectutem aegritudine animi affigere vereor; alias
crabrones monadicos per totam Germaniam irritare non perhorrescerem.
}
\end{itemize}

\begin{itemize}
\item[]{\small
I confess that I avoided all this also by the reason
that I did not want to look as a bragger attacking the writhing of the great scholars 
rather the a man pursuing  truth.
The same reason has been preventing for a long time
to submit for consideration of the scientific council my views of monads.
Although I am absolutely convinced that this mystical teaching must
be completely destroyed  by my arguments, I am afraid to spoil the elder
years of the scholar  whose benefactions to me I could not forget;
otherwise I would not be scared  to tease hornet-monadists throughout the whole
of Germany.
}
\end{itemize}

It is worth observing that Lomonosov had in mind not the idea of Leibniz himself but rather
the exposition of monadism in the writings of Wolff and his numerous descendants.
This day we know  Wolff's letter to Ernst Christoph von Manteuffel as of May 11, 1746
which shows that Wolff considered his metaphysics as different from that of Leibniz 
(cp.~\cite{Lenders}).

Let us make a mental ``physicalistic'' experiment
and aim  a strong microscope at
a region about  a point at a mathematical line.
We will see in the
eyepiece a~blurred and dispersed cloud with unclear frontiers
which is a visualization of the point under investigation.
Under greater magnification,
the portion of the ``point-monad'' we are looking at
will enlarge, revealing extra details whereas
disappearing  partially from sight.
However, we are still inspecting the same standard real number
which  you might prefer to percept as
described by this process of
``studying the microstructure of a~physical straight line.''
Visualizing a point by microscope reveals
its monadic essence. Lomonosov and even Leibniz
could reason  likewise or approximately so.
In any case, the view of the monad of a standard real number
as the collection of all infinitely close points
is generally adopted in the contemporary infinitesimal
analysis resurrected under the name of {\it nonstandard analysis}
in the works by Abraham Robinson in 1961 (cp.~\cite{Rob8}, \cite{GorKusKut})..

More than two dozens of decades elapsed from the death of
Mikha\u\i{}lo Lomonosov, but his creative contribution still
inspires thought in connection with the most topical and brand-new areas of mathematics and
natural sciences.
His enviable fate gives a supreme example for drafting life.

%\vfill\eject

\bibliographystyle{plain}

\end{document}